\pdfoutput=1

\documentclass[a4paper,10.5pt]{article}
\usepackage[latin1]{inputenc}
\usepackage{amsmath}
\usepackage{amsfonts}
\usepackage{amssymb}
\usepackage{geometry}
\usepackage{graphicx}

\newtheorem{theo}{Th\'eor\`eme}[section]
\newtheorem{cor}[theo]{Corollaire}
\newtheorem{prop}[theo]{Proposition}
\newtheorem{lemme}[theo]{Lemme}
\newtheorem{pte}[theo]{Propriétés}
\newtheorem{defi}[theo]{D\'efinition}
\newtheorem{rema}[theo]{Remarque}
\newtheorem{exe}[theo]{Exemple}
\newtheorem{conj}[theo]{Conjecture}

\newcommand{\cA}{\mathcal{A}}
\newcommand{\cM}{\mathcal{M}}
\newcommand{\cS}{\mathcal{S}}
\newcommand{\cR}{\mathcal{R}}
\newcommand{\cU}{\mathcal{U}}
\newcommand{\bcA}{\boldsymbol{\cA}}
\newcommand{\bcS}{\boldsymbol{\cS}}

\newcommand{\wcM}{\widetilde\cM}
\newcommand{\wcU}{\widetilde\cU}

\newcommand{\VT}{\Vert T_n \Vert}
\newcommand{\VcM}{\Vert \cM_n \Vert}
\newcommand{\VwcM}{\Vert \widetilde\cM_n \Vert}

\newcommand{\bpi}{\boldsymbol{\pi}}
\newcommand{\bmu}{\boldsymbol{\mu}}
\newcommand{\N}{\mathbb{N}}
\newcommand{\Z}{\mathbb{Z}}
\newcommand{\ov}{\overline}
\newcommand{\wh}{\widehat}
\newcommand{\sn}{\sqrt{n}}
\newcommand{\pr}{\textit{Preuve : }}

\begin{document}

\author{Jean-Paul Cardinal}

\title{Une suite de matrices symétriques en rapport avec la fonction de Mertens}
\maketitle
\begin{abstract} 
In this paper we explore a class of equivalence relations over $\N^\ast$ from which is
constructed a sequence of symetric matrices related to the Mertens function. 
From numerical experimentations we formulate a conjecture, about 
the growth of the quadratic norm of these matrices, which implies the Riemann conjecture.
This suggests that matrix analysis methods may play a more important part in this classical and difficult problem.
\end{abstract}

\section{Introduction.}

Rappelons les définitions, sur les entiers naturels,
 de la fonction de Möbius :
 \begin{equation}
\mu (k)=\left\{
\begin{array}{cc}
(-1)^r   & \text{ si k est un produit de r facteurs premiers distincts,}\\ 
0  & \text{ si k  est  divisible  par  un  carré.}
\end{array}
\right.
\end{equation}
et de la fonction de Mertens :
 \begin{equation}
\label{eq:mertens}
M(n)=\sum_{1\le k\le n} \mu(k)
\end{equation}

Dans cet article, nous construisons une suite de matrices carrées symétriques 
$\cM=\left( \cM_n\right)_{n\in\N^\ast} $ vérifiant pour tout $n$ :
\begin{equation}
\vert M(n) \vert \le \VcM
\end{equation}
où dans toute la suite, $\Vert A \Vert$ désigne $\Vert A \Vert_2$ la norme quadratique sur les matrices, (appelée aussi $2$-norme, voir par exemple \cite{GL} p.56).
L'importance d'une majoration de la suite $\vert M(n) \vert$ tient au fait
 que la conjecture
 \begin{equation}
\forall \epsilon > 0, M(n)=O(n^{1/2+\epsilon})
\end{equation}
est équivalente à l'hypothèse de Riemann (cf. Littlewood \cite{L}). 
F. Mertens \cite{M} conjectura dès $1897$ que $\vert M(n)\vert \le n^{1/2}$ pour tout $n$. Stieljes annonça en $1885$ avoir montré que $M(n) n^{-1/2}$ est bornée, mais sans publier de preuve.
La conjecture de Mertens fut réfutée en $1985$ par H.M. Odlyzko et H.J.J. te Riele \cite{OR}, mais leur preuve est non-constructive, et on ne connait toujours pas aujourd'hui de valeur $n$ explicite vérifiant $\vert M(n)\vert > n^{1/2}$, bien qu'on ait montré qu'une telle valeur $n$ existe
 dans l'intervalle $\left[ 10^{14} , 10^{1.4\times 10^{64}} \right] $, (cf. \cite{KL1} et \cite{P}). \\
 
 La construction de la suite de matrices $\cM$ est détaillée dans la section $2$.
 Dans la section $3$ nous examinons la croissance en norme de deux autres suites de matrices présentant un lien avec $\cM$.
 Dans la section $4$ nous effectuons, pour un ensemble de valeurs de $n$, l'évaluation numérique en Octave de $\VcM$, dont les résultats expérimentaux nous incitent à conjecturer que :
\begin{conj} 
\label{maj_fondamentale}
\begin{equation}
 \forall \epsilon > 0, \VcM =O(n^{1/2+\epsilon})
\end{equation}
\end{conj}

\section{Construction de la suite $\cM$.}

\subsection{Relation d'équivalence sur $\N^\ast$; monoïde quotient $\wh\cS=\N^\ast/\cR$.}

\begin{defi}
Pour $n\in \N^\ast$ fixé, nous définissons sur $\N^\ast$ la relation d'équivalence :
$$i \ \cR \ j \Longleftrightarrow \left[ n/i \right]=\left[ n/j \right]$$
\end{defi}

\begin{exe}
Pour $n=16$ nous obtenons les huit classes d'équivalence :\\
$\lbrace 1\rbrace,\lbrace 2\rbrace,\lbrace 3\rbrace,\lbrace 4 \rbrace,\lbrace 5\rbrace ,
\lbrace 6,7,8\rbrace,\lbrace 9,10,11,12,13,14,15,16 \rbrace$ et $\lbrace 17,18,\cdots\rbrace$.
\end{exe}

Considérant la structure en intervalles des classes d'équivalences, nous pouvons identifier sans ambiguité 
chaque classe par son plus grand représentant, en convenant de noter  $\infty$ le plus grand représentant de la classe non bornée. Pour bien les distinguer nous typographions les représentants en caractères simples et les classes d'équivalence en caractères surmontés d'un chapeau. Notons $\ov\cS$ l'ensemble de ces plus grands représentants et $\cS=\ov\cS\setminus\lbrace \infty \rbrace$. Notons aussi $\wh\cS$ l'ensemble des classes, c'est-à-dire $\wh\cS=\N^\ast/\cR$. 

\begin{exe}
\label{exeS}
 Pour n = 16 :\\
$\ov\cS=\{1,2,3,4,5,8,16,\infty \}$, 
$\cS=\{1,2,3,4,5,8,16 \}$,
$\wh\cS=\{\wh{1},\wh{2},\wh{3},\wh{4},\wh{5},\wh{8},\wh{16},\wh{\infty} \}$, avec 
${\wh 1}=\lbrace 1\rbrace,\\{\wh 2}=\lbrace 2\rbrace,
{\wh 3}=\lbrace 3\rbrace,{\wh 4}=\lbrace 4 \rbrace,{\wh 5}=\lbrace 5\rbrace ,{\wh 8}=\lbrace 6,7,8\rbrace, {\wh 16}=\lbrace 9,10,11,12,13,14,15,16 \rbrace$, $ \wh{\infty}=\lbrace 17,18,\cdots \rbrace$.
\end{exe}

\begin{prop}
\label{diezS}
Soit $n\in \N^\ast$ et $\cS$ défini comme ci-dessus, c'est-à-dire l'ensemble des plus grands éléments des classes finies de la relation d'équivalence $\cR$. Par ailleurs, pour chaque $k\in\N^\ast$, posons $\ov{k}=[n/k]$.
\begin{enumerate}
\item On a l'inclusion :
$\left\lbrace 1,2,\cdots,\left[ \sn\right] \right\rbrace \subset \cS$.
\item Pour tout $k$ dans $\cS$, on a $\ov{k}\in\cS$ et $\ov{\ov{k}}=k$, c'est-à-dire que l'application 
$k \mapsto \ov{k}$ est une involution décroissante de $\cS$. En fait $k \mapsto \ov{k}$ est simplement l'application qui renverse l'ordre des éléments de $\cS$.
\item
Il existe une constante $C$ telle que pour tous $n\in\N^\ast$ et $k\in\cS$, on a $\dfrac{k^+}{k}\le C$,
où $k^+$ désigne le successeur de $k$ dans $\cS$.
\item On a l'estimation suivante de $\#\cS$ en fonction de $n$ :
\begin{equation}
\begin{array}{c}
\text{Si } n<[\sn]^2+[\sn] \text{ alors } \cS=\left\lbrace 1,\cdots,[\sn]=\ov{[\sn]},\cdots,n=\ov{1} \right\rbrace {\mathrm\ et\  } \#\cS=2[\sn]-1 \\
\text{Si } n\ge [\sn]^2+[\sn] \text{ alors } \cS=\left\lbrace 1,\cdots,[\sn],\ov{[\sn]},\cdots,n=\ov{1} \right\rbrace  {\mathrm\ et\  } \#\cS=2[\sn]
\end{array}
\end{equation}
 \end{enumerate}
\end{prop}

\pr
\begin{enumerate}
\item
Nous commençons par montrer que chaque singleton $\left\lbrace k\right\rbrace $, avec 
$1\le k <\sn$, est une classe d'équivalence. En effet, si $k<[\sn ]$ et donc $k+1\le [\sn ]$,
alors $\dfrac{n}{k}-\dfrac{n}{k+1}=\dfrac{n}{k(k+1)}>1$ ce qui entraine que 
$\left[ \dfrac{n}{k+1}\right] <\left[ \dfrac{n}{k}\right] $,
 c'est-à-dire que $k$ et $k+1$ ne sont pas dans la même classe.\\
Traitons maintenant le cas $k=\left[ \sn\right] $, en distinguant deux cas :
\begin{itemize}
\item si $\dfrac{n}{k}-\dfrac{n}{k+1}=\dfrac{n}{k(k+1)} \ge 1$ alors
$\left[ \dfrac{n}{k+1}\right] <\left[ \dfrac{n}{k}\right] $.
\item si $\dfrac{n}{k}-\dfrac{n}{k+1}=\dfrac{n}{k(k+1)} < 1$ alors 
$\displaystyle{\frac{n}{k+1}<k\le\frac{n}{k}}$ et comme $k$ est entier,
 cela entraine bien que $\left[ \dfrac{n}{k+1}\right] <\left[ \dfrac{n}{k}\right] $.
\end{itemize}
Dans les deux cas, on en déduit que $[\sn]$ et $[\sn]+1$ ne sont pas dans la même classe d'équivalence.
\item
Partons de $\ov{k}=[n/k]$, c'est-à-dire $k\ov{k}\le n < k\ov{k}+k$ ce qui se réécrit :
$\dfrac{n}{\ov{k}+1} < k \le \dfrac{n}{\ov{k}}$.\\
$k$ étant entier, ceci entraine que 
$\displaystyle{\left[ \frac{n}{\ov{k}+1}\right] < \left[ \frac{n}{\ov{k}}\right]}$
et donc que $\ov{k}\in\cS$.\\
Soit $k\in\cS$. On a $\left[ \dfrac{n}{k+1}\right] <\left[ \dfrac{n}{k}\right] $ 
ce qui entraine que
$ \dfrac{n}{k+1} <\left[ \dfrac{n}{k}\right] =\ov{k}$ c'est-à-dire
$n<k\ov{k}+\ov{k}$. Par ailleurs on a $k\ov{k}\le n$ d'où la double inégalité
$\displaystyle{k\le \frac{n}{\ov{k}} <k+1}$ qui signifie que $k=[n/\ov{k}]=\ov{\ov{k}}$.
\item Distinguons deux cas :
\begin{itemize}
\item
Si $k\le [\sn]$ alors $\dfrac{k^+}{k}=\dfrac{k+1}{k}\le 2$. 
\item
Si $k > [\sn]$, posons $l=\ov{k}$ et notons $l^-$ le prédécesseur de $l$ dans $\cS$.\\
 Alors
$\dfrac{k^+}{k}=\dfrac{\ov{l^-}}{\ov{l  }}\le \dfrac{[n/l^-]}{[n/l]}\le \dfrac{l}{l^-}\dfrac{n}{n-l}$.
Puisque $l\le [\sn]$ on peut comme précédemment majorer $\dfrac{l}{l^-}$ par $2$ et puisque la quantité
$\dfrac{n}{n-\sn}$ est une fonction décroissante de $n$ on peut la majorer par 
$\dfrac{2}{2-\sqrt{2}}=2+\sqrt{2}$.
\end{itemize}
On a donc la conclusion recherchée avec $C=4+2\sqrt{2}$. On pourrait d'ailleurs s'apercevoir pratiquement que $C=3$ convient.
\item 
Supposons $k > [\sn]$. On a alors successivement $k \ge \sn$, $n/k \le \sn$ et $\ov{k} \le [\sn]$.
En d'autres termes, $[\sn]$ est le plus grand élément $k$ de $\cS$ tel que $k\le \ov{k}$. En utilisant le fait que $k \mapsto \ov{k}$ est une involution décroissante sur $\cS$, et en distinguant les cas
$[\sn] = \ov{[\sn]}$ et $[\sn] < \ov{[\sn]}$, on en déduit la forme escomptée de $\cS$.\\
Pour finir cherchons à réécrire la condition $[\sn] = \ov{[\sn]}$. Posons $\sn=k+\alpha$ avec $k$ entier et 
$0\le \alpha < 1$.
On a :
$$n = k^2 + 2\alpha k +\alpha^2 \Leftrightarrow n/k = k + 2\alpha +\alpha^2/k$$ donc :
$$[\sn] = \ov{[\sn]} \Leftrightarrow 2\alpha +\alpha^2/k < 1 \Leftrightarrow \alpha^2+2k\alpha-k <0
\Leftrightarrow n < k^2 +k$$
 \end{enumerate}
\hfill$\square$

\begin{rema}
On peut synthétiser les deux cas du point 4. de la proposition ci-dessus sous la forme suivante, valable pour tout $n\in\N^\ast$ :
\begin{equation}
\#\cS=[\sn]+[\sqrt{n+1/4}-1/2]
\end{equation}
\end{rema}

\begin{lemme}
On a d'une façon générale, pour des entiers non nuls $i,j,n$ :
\begin{equation}
\left[\left[ n/i \right]/j\right]=\left[ n/ij \right]
\end{equation}
\end{lemme}
\pr Posons $n/i=u+\alpha$, où $u$ est entier et $0\le \alpha <1$. On a donc 
$[n/i]/j=u/j$ et $n/ij=u/j+\alpha/j$ ce qui entraine l'inégalité 
$\left[\left[ n/i \right]/j\right]\le \left[ n/ij \right]$.
 Si cette dernière inégalité était stricte cela entrainerait l'existence d'un entier $v$ tel que
  $u/j < v \le u/j+\alpha/j $
donc $u < vj \le u+\alpha < u+1$, ce qui est impossible étant donné que $u$ et $vj$ sont entiers.
\hfill$\square$

\begin{prop}
\label{multiR}
 $\cR$ est compatible avec la multiplication sur $\N^\ast$, au sens où pour tous $i,j,k \in \N^\ast$  on a $i \ \cR \ j \Longrightarrow ik \ \cR \ jk$. On peut donc définir sur $\wh \cS$ une multiplication induite, par la formule $\wh i\wh j=\wh{ij}$, où rappelons le, $\wh i$ désigne la classe d'équivalence de l'entier $i$ et $\wh \infty$ la classe de tout entier strictement supérieur à $n$.
\end{prop} 
\pr
On suppose $i \ \cR \ j$ donc $\left[ n/i \right]=\left[ n/j \right]$. 
Alors d'après le lemme précédent 
$[ n/ik ]=\left[ [n/i]/k \right] =\left[ [n/j]/k \right] =[ n/jk]$,
 c'est-à-dire $ik \ \cR \ jk $.
\hfill$\square$
\\

$\N^\ast$, muni de la multiplication usuelle, est un monoïde
 (ensemble muni d'une loi interne associative et disposant d'un élément neutre) commutatif. 
  $\cR$ étant compatible avec la multiplication sur $\N^\ast$,
 l'ensemble quotient $\wh\cS=\N^\ast/\cR$, muni de la multiplication induite, est donc lui aussi un monoïde commutatif.

\begin{exe}
Table de multiplication du monoïde $\wh\cS=\N^\ast/\cR$ pour $n=16$.
\begin{center}
\begin{footnotesize}
$\begin{array}{c|cccccccc}
    &   \wh{1} &   \wh{2} &   \wh{3} &   4 &   \wh{5} &   \wh{8} &   \wh{16} &  \wh{\infty} \\ \hline
  \wh{1} &   \wh{1} &   \wh{2} &   \wh{3} &   \wh{4} &   \wh{5} &   \wh{8} &   \wh{16} & \wh{\infty} \\ 
  \wh{2} &   \wh{2} &   \wh{4} &   \wh{8} &   \wh{8} &   \wh{16} &   \wh{16} &  \wh{\infty} &  \wh{\infty} \\ 
  \wh{3} &   \wh{3} &   \wh{8} &   \wh{16} &   \wh{16} &   \wh{16} &  \wh{\infty} &   \wh{\infty} &  \wh{\infty}\\
  \wh{4} &   \wh{4} &   \wh{8} &   \wh{16} &   \wh{16} &  \wh{\infty} &  \wh{\infty} &  \wh{\infty} &  \wh{\infty}\\
  \wh{5} &   \wh{5} &   \wh{16} &   \wh{16} &  \wh{\infty} &  \wh{\infty} &  \wh{\infty} &  \wh{\infty} &  \wh{\infty} \\
  \wh{8} &   \wh{8} &   \wh{16} &  \wh{\infty} &  \wh{\infty} &  \wh{\infty} &  \wh{\infty} &  \wh{\infty} &  \wh{\infty} \\
  \wh{16} &   \wh{16} &  \wh{\infty} &  \wh{\infty} &  \wh{\infty} &  \wh{\infty} &  \wh{\infty} &  \wh{\infty} &  \wh{\infty} \\
\wh{\infty} & \wh{\infty} &  \wh{\infty} &  \wh{\infty} &  \wh{\infty} &  \wh{\infty} &  \wh{\infty} &  \wh{\infty} &  \wh{\infty} 
\end{array}$
\end{footnotesize}
\end{center}
\end{exe}

\subsection{Algèbres $A=\Z^{\N^\ast}$, $\wh{A}=\Z^{\wh\cS}$; algèbre quotient $\bcA$.}
Lorsqu'on dispose d'un monoïde, on peut construire sa $\Z$-algèbre de monoïde en formant toutes les combinaisons linéaires formelles, finies (ou infinies lorsque les produits sont naturellement bien définis), à coefficients entiers relatifs, d'éléments du monoïde.
 La multiplication dans l'algèbre se fait d'une manière évidente, en utilisant la multiplication dans le monoïde et la distributivité de la multiplication sur l'addition.

\subsubsection{L'algèbre $A$ du monoïde $\N^\ast$.}
Si l'on choisit le monoïde $\N^\ast$ muni de la multiplication usuelle, alors l'algèbre ${A}=\Z^{\N^\ast}$ est l'algèbre des séries de Dirichlet, ici à coefficients entiers, (voir par exemple \cite{KM} p. 111). Le produit $\star$ sur $A$ s'appelle produit de convolution ou produit de Dirichlet. Voici quelques propriétés bien connues de $A$ :

\begin{pte}
\label{dirichlet}
\begin{enumerate}
\item
Si ${ a}=(a_1,\cdots,a_k,\cdots)$, ${ b}=(b_1,\cdots,b_k,\cdots)$ sont des éléments de ${A}$,\\
 et ${ c}=(c_1,\cdots,c_k,\cdots)=a\star b$, alors $c_k=\sum_{ij= k} a_i b_j$. Si $e_i,e_j$ sont les $i$-ème et $j$-ème vecteurs de la base canonique de ${A}$, alors $e_i\star e_j= e_{ij}$.
\item
L'élément unité est $e_1=(1,0,\cdots,0,\cdots)$. Un élément ${ a}=(a_1,\cdots,a_n,\cdots)$ est inversible si et seulement si $a_1 = \pm 1$. 
\item
Si $u=(1,1,\cdots,1,\cdots)$, alors $u^{-1}=\mu$, où $\mu$ est la suite de Möbius, (cf. \cite{KM} p. 112).
\end{enumerate}
\end{pte}

\subsubsection{L'algèbre $\wh A$ du monoïde $\wh\cS$, l'algèbre quotient $\bcA$.}
Si l'on choisit le monoïde $\wh\cS = \N^\ast/\cR$ muni de la multiplication induite, alors l'algèbre ${\wh A}=\Z^{\wh\cS}$ est une $\Z$-algèbre de dimension $\#\wh\cS$, dont une base est $\wh\cS$.
Il s'avère cependant plus intéressant pour nous de considérer l'algèbre quotient
$\bcA = \wh A /{\wh\infty{\wh A}}$ où  $\wh\infty{\wh A} = \Z\wh\infty$ est l'idéal principal de $\wh A$ engendré par $\wh\infty$.
Soit $\varpi$ la projection canonique de ${\wh A}$ sur $\bcA$. Nous noterons les images par $\varpi$ des éléments de la base $\wh\cS$, et d'une façon générale les vecteurs de $\bcA$, par des lettres en gras. Ainsi par exemple 
$\varpi(\wh k)={\bf k}$, $\varpi(\wh\infty)={\bf 0}$.
Lorsque $k$ parcourt l'ensemble $\cS$, (cf. exemple \ref{exeS} et proposition \ref{diezS}), $\bf k$ parcourt un ensemble que nous noterons $\bcS$ et qui est une base de $\bcA$. $\bcS$, qui est l'image par $\varpi$ de la base $\wh\cS$, sera dorénavant appelée base canonique de $\bcA$. On a bien sûr $\#\bcS=\#\cS$ quantité qui a été calculée en fonction de $n$ dans la proposition \ref{diezS}.
Avec ces notations enfin, il est facile de construire la table de multiplication de la base $\bcS$ de l'algèbre $\bcA$, à partir de la table de multiplication de $\wh \cS$. Voici comment :
\begin{itemize}
\item
supprimer la dernière ligne et la dernière colonne de la table de $\wh \cS$ et remplacer les $\wh\infty$ restants par des ${\bf 0}$ (ceci traduit le fait que $\varpi(\wh\infty)={\bf 0}$).
\item
supprimer les chapeaux et réécrire les entiers en caractères gras (ceci correspond à la réécriture $\varpi(\wh k)={\bf k}$ lors de la projection sur $\bcA$).
\end{itemize}

\pagebreak
\begin{exe}
Table de multiplication de la base $\bcS$ de $\bcA$ pour $n=16$.\\
\begin{center}
\begin{footnotesize}
\boldmath
$\begin{array}{c|ccccccc}
    &   1 &   2 &   3 &   4 &   5 &   8 &   16  \\ \hline
  1 &   1 &   2 &   3 &   4 &   5 &   8 &   16  \\ 
  2 &   2 &   4 &   8 &   8 &   16 &   16 &  0  \\ 
  3 &   3 &   8 &   16 &   16 &   16 &  0 &   0 \\
  4 &   4 &   8 &   16 &   16 &  0 &  0 &  0 \\
  5 &   5 &   16 &   16 &  0 &  0 &  0 &  0  \\
  8 &   8 &   16 &  0 &  0 &  0 &  0 &  0  \\
  16 &   16 &  0 &  0 &  0 &  0 &  0 &  0  
\end{array}$
\end{footnotesize}
\end{center}
\end{exe}

\subsection{Morphisme naturel $\bpi$ de $\Z^{\N^\ast}$ sur ${\bcA}$.}

\begin{prop}
$n$ étant toujours fixé, et $\cR$ la relation d'équivalence associée, rappelons que $e_i$ désigne le $i$ème vecteur canonique de l'algèbre $\Z^{\N^\ast}$ et définissons le morphisme de $\Z$-modules :
\begin{equation}
\vartheta :
\left\vert
\begin{array}{c}
\Z^{\N^\ast} \mapsto \wh A \\
e_i \mapsto \wh i
\end{array} 
\right. 
\end{equation}
Alors $\vartheta$ est un morphisme de $\Z$-algèbres, c'est-à-dire que pour tous $a,b$ dans $\Z^{\N^\ast}$, on a :
\begin{equation}
\vartheta(a\star b)=\vartheta(a)\vartheta(b)
\end{equation}
\end{prop}
\pr
Par linéarité, il suffit de le montrer pour deux éléments $e_i,e_j$ de la base canonique
de $\Z^{\N^\ast}$. D'après les propositions \ref{multiR} et \ref{dirichlet} on a :
$$\vartheta(e_i\star e_j)=\vartheta(e_{ij})=\wh{ij}=\wh i\wh j=\vartheta(e_i)\vartheta(e_j)$$
\hfill$\square$

\begin{cor}
\label{bpi_morph}
Le composé $\bpi=\varpi \circ \vartheta$ :
\begin{equation}
\bpi :\left\vert  \Z^{\N^\ast} \rightarrow \bcA \right. 
\end{equation}
est un morphisme de $\Z$-algèbres.
\end{cor}

\begin{prop}
L'image de $a=(a_1,\cdots,a_\kappa,\cdots)\in\Z^{\N^\ast}$ par $\bpi$ vaut :
\begin{equation}
\bpi(a) = \sum_{ k\in \cS}\left( \sum_{\kappa\in {\wh k}}a_\kappa\right) {\bf k}
\end{equation}
\end{prop}
\pr
Cette proposition est évidente dès que l'on a remarqué que pour chaque $k\in \cS$, 
 $\wh k$ est une partie finie de $\N^\ast$ et que par suite la somme $\sum_{\kappa\in {\wh k}}a_\kappa$ est bien définie.
\hfill$\square$

\begin{cor}
\label{mertens}
Par le morphisme $\bpi$, 
\begin{itemize}
\item
 $u=(1,1,\cdots,1,\cdots)$ est transformé en ${\bf u}=(\#{\wh k})_{k\in \cS}$.
  Si pour chaque $k\in \cS$
 on note $k^-$ le prédécesseur de $k$ dans $\cS$, avec la convention $1^-=0$, alors
on a  $\#{\wh k}=k-k^-$ et donc :
\begin{equation}
{\bf u}=(k-k^-)_{k\in\cS}
\end{equation}
\item
 $\mu=(\mu(1),\mu(2),\cdots,\mu(k),\cdots)$ est transformé en 
${\bmu}=\left(\sum_{\kappa\in {\wh k}}\mu(\kappa)\right) _{ k\in \cS}$,
c'est-à-dire :
\begin{equation}
\bmu=\left( M(k)-M(k^-)\right) _{k\in\cS}
\end{equation}
où $M$ désigne la fonction de Mertens sur les entiers.
\end{itemize}
\end{cor}

\begin{exe}
\label{mertens_exe}
Pour $n=16$ on a  :
\begin{itemize}
\item
${\bf u}=(1,1,1,1,1,3,8)={\bf 1} + {\bf 2} + {\bf 3} + {\bf 4} + {\bf 5} + 3\ {\bf 8} + 8\ {\bf 16}$.
\item
${\bmu}=(1,-1,-1,0,-1,0,1)={\bf 1-2-3-5+16}\\ 
=\mu(1){\bf 1}+\mu(2) {\bf 2}+\mu(3) {\bf 3}+\mu(4) {\bf 4}+\mu(5) {\bf 5}
+\left( M(8)-M(5)\right)  {\bf 8}+\left( M(16)-M(8)\right) {\bf 16}$.
\end{itemize}
\end{exe}

\subsection{Représentation linéaire de l'algèbre $\bcA$.}
Pour chaque ${\bf a} \in \bcA$, l'application :
\begin{equation}
\label{eq:linear_rep}
\rho :
\left\vert
\begin{array}{c}
\bcA \mapsto \bcA \\ 
{\bf x} \mapsto {\bf a  x}
\end{array} 
\right.
\end{equation}
est linéaire, et admet donc dans $\bcS$, la base canonique de $\bcA$, une 
matrice $\rho({\bf a})$. Notons $s=\#\bcS=\#\cS$ la dimension de $\bcA$ ( $s=\#\cS$ a été calculé en fonction de $n$ dans la proposition $\ref{diezS}$).
Si $\cM_s(\Z)$ désigne l'algèbre des matrices carrées de taille $s$ à coefficients dans $\Z$, alors l'application
\begin{equation}
\label{eq:linear_iso}
\rho:
\left\vert
\begin{array}{c}
{\bcA}\mapsto \mathcal{M}_s(\Z) \\ 
{\bf a} \mapsto \rho({\bf a})
\end{array} 
\right.
\end{equation}
est ce qu'on appelle la représentation régulière de l'algèbre $\bcA$, représentation qui est fidèle, c'est-à-dire que le morphisme $\rho$ est injectif.
L'ensemble des matrices $\rho({\bf a}), {\bf a} \in \bcA$ est donc une sous-algèbre commutative de dimension $s$ de $\cM_s(\Z)$, dont une base est formée
des matrices $\rho({\bf k}),\ {\bf k}\in\bcS$. Enfin, puisqu'il y a une bijection très naturelle entre $\cS$ et $\bcS$, nous choisirons $\cS$ comme ensemble 
d'indexation pour les lignes et les colonnes des matrices $\rho({\bf a})$. Par exemple pour $n=16$
la dernière colonne d'une matrice $\rho({\bf a})$ ne sera pas la colonne d'indice $7$, mais la colonne 
d'indice $16$.

\begin{exe}
\label{representants}
 Pour n=16, outre $\rho({\bf 1})$ qui est l'identité, nous obtenons les matrices représentant \\
 ${\bf 2, 3, 4, 5, 6, 8, 16}$ (où la plupart des entrées nulles sont laissées en blanc pour plus de lisibilité) :
\begin{center}
\begin{footnotesize}
$
\begin{array}{c|ccccccc|}
\rho({\bf 2}) & { 1} & { 2} & { 3} & { 4} & { 5} & { 8} & { {16}}\\ \hline
{ 1} & & & & & & &\\ 
{ 2} & 1 & & & & & &\\ 
{ 3} & & & & & & &\\ 
{ 4} & & 1 & & & & &\\ 
{ 5} & & & & & & &\\ 
{ 8} & & & 1 & 1 & & &\\ 
{ {16}} & & & & & 1 & 1 &\\ \hline
\end{array} 
$
$
\begin{array}{c|ccccccc|}
\rho({\bf 3}) & { 1} & { 2} & { 3} & { 4} & { 5} & { 8} & { {16}}\\ \hline
{ 1} & & & & & & &\\ 
{ 2} & & & & & & &\\ 
{ 3} & 1 & & & & & &\\ 
{ 4} & & & & & & &\\ 
{ 5} & & & & & & &\\ 
{ 8} & & 1 & & & & &\\ 
{ {16}} & & & 1 & 1 & 1 & &\\ \hline
\end{array} 
$
$
\begin{array}{c|ccccccc|}
\rho({\bf 4}) & { 1} & { 2} & { 3} & { 4} & { 5} & { 8} & { {16}}\\ \hline
{ 1} & & & & & & &\\ 
{ 2} & & & & & & &\\ 
{ 3} & & & & & & &\\ 
{ 4} & 1 & & & & & &\\ 
{ 5} & & & & & & &\\ 
{ 8} & & 1 & & & & &\\ 
{ {16}} & & & 1 & 1 & & &\\ \hline
\end{array} 
$\\
$
\begin{array}{c|ccccccc|}
\rho({\bf 5}) & { 1} & { 2} & { 3} & { 4} & { 5} & { 8} & { {16}}\\ \hline
{ 1} & & & & & & &\\ 
{ 2} & & & & & & &\\ 
{ 3} & & & & & & &\\ 
{ 4} & & & & & & &\\ 
{ 5} & 1 & & & & & &\\ 
{ 8} & & & & & & &\\ 
{ {16}} & & 1 & 1 & & & &\\ \hline
\end{array} 
$
$
\begin{array}{c|ccccccc|}
\rho({\bf 8}) & { 1} & { 2} & { 3} & { 4} & { 5} & { 8} & { {16}}\\ \hline
{ 1} & & & & & & &\\ 
{ 2} & & & & & & &\\ 
{ 3} & & & & & & &\\ 
{ 4} & & & & & & &\\ 
{ 5} & & & & & & &\\ 
{ 8} & 1& & & & & &\\ 
{ {16}} & & 1& & & & &\\ \hline
\end{array} 
$
$
\begin{array}{c|ccccccc|}
\rho({\bf {16}}) & { 1} & { 2} & { 3} & { 4} & { 5} & { 8} & { {16}}\\ \hline
{ 1} & & & & & & &\\ 
{ 2} & & & & & & &\\ 
{ 3} & & & & & & &\\ 
{ 4} & & & & & & &\\ 
{ 5} & & & & & & &\\ 
{ 8} & & & & & & &\\ 
{ {16}} & 1& & & & & &\\ \hline
\end{array} 
$
\end{footnotesize}
\end{center}
et suivant l'exemple \ref{mertens_exe}, les matrices $\rho({\bf u})$ et $\rho({\bmu})$ :
\begin{center}
\begin{footnotesize}
$
\begin{array}{c|ccccccc|}
\rho({\bf u}) & { 1} & { 2} & { 3} & { 4} & { 5} & { 8} & { {16}}\\ \hline
{ 1} & 1& & & & & &\\ 
{ 2} & 1& 1& & & & &\\ 
{ 3} & 1& 0& 1& & & &\\ 
{ 4} & 1& 1& 0& 1& & &\\ 
{ 5} & 1& 0& 0& 0& 1& &\\ 
{ 8} & 3& 2& 1& 1& 0& 1&\\ 
{ {16}} & 8& 4& 3& 2& 2& 1& 1\\ \hline
\end{array} 
$
\hspace{1cm}
$
\begin{array}{c|ccccccc|}
\rho({\bmu}) & { 1} & { 2} & { 3} & { 4} & { 5} & { 8} & { {16}}\\ \hline
{ 1} & 1  &   &    &   &   &   &\\ 
{ 2} & -1 & 1 &    &   &   &   &\\ 
{ 3} & -1 &  0 &  1 &   &   &   &\\ 
{ 4} &  0  &-1 &  0  & 1 &   &   &\\ 
{ 5} & -1 &  0 &  0  &  0 & 1 &   &\\ 
{ 8} &  0  & -1& -1 & -1& 0  & 1 &\\ 
{ {16}}& 1  & -1& -2 & -1& -2& -1& 1\\ \hline
\end{array} 
$
\end{footnotesize}
\end{center}
sur lesquelles on constate bien que pour ${\bf a}\in {\bcA}$, les coefficients de ${\bf a}$ dans la base $\bcS$ se retrouvent dans la première colonne de $\rho({\bf a})$, (cf. exemple \ref{mertens_exe}).
\end{exe}

\begin{defi}
Pour $n$ fixé, soit $T$ la matrice de taille $s=\#\cS$, symétrique, qui vaut $1$ au dessus de la perdiagonale, et $0$ strictement en dessous de la perdiagonale.
\end{defi}

\begin{exe} Pour $n=16$ :
$$T=
\begin{small}
\left[
\begin{array}{ccccccc}
  1& 1& 1& 1& 1& 1 &1\\ 
  1& 1& 1& 1& 1& 1 &\\ 
  1& 1& 1& 1& 1&   &\\ 
  1& 1& 1& 1&  &   &\\ 
  1& 1& 1&  &  &   &\\ 
  1& 1&  &  &  &   &\\ 
  1&  &  &  &  &   &\\
\end{array} 
\right]
\end{small}
$$
\end{exe}

\begin{prop}
\label{TM_sym} 
Pour tout ${\bf a}\in {\bcA}$, $T\rho({\bf a})$ est symétrique.
\end{prop}
\pr
Il suffit de le démontrer si $\bf a$ est égal à l'un des vecteurs de base $\bf k\in \bcS$.
Soit donc ${\bf k}\in \bcS$. La colonne $v$ d'indice ${ j}$ de la matrice $\rho(\bf k)$
contient, outre des $0$, un unique $1$ à l'indice ${ l}$ tel que ${\bf l = kj}$
c'est-à-dire que l'on a $[n/l]=[n/jk]$.
La colonne $Tv$ est donc la colonne d'indice ${ l}$ de $T$,
c'est-à-dire constituée de $1$ pour tous les indices ${ i}$ tels que $i\le \ov{l}$
et de $0$ au delà, où $\ov{l}=[n/l]$ est le symétrique de $l$ dans la liste $\cS$, (proposition \ref{diezS}). On en déduit que, $v_i$ désignant la composante d'indice $i$ de $v$ :
$$
v_{ i}=1 \Leftrightarrow i\le [n/l]=[n/jk]
\Leftrightarrow i\le n/jk
\Leftrightarrow  ij\le n/ k
\Leftrightarrow  ij\le [n/ k]$$
condition qui est bien symétrique en les indices $i,j$, ce qu'il fallait démontrer.
\hfill$\square$

\begin{exe}
\label{TuTmu_exe}
Pour $n=16$, on a 
$T\rho({\bf 1})=T$, et 

\begin{center}
\begin{footnotesize}
$
\begin{array}{c|ccccccc|}
{T\rho(\bf 2}) & { 1} & { 2} & { 3} & { 4} & { 5} & { 8} & { {16}}\\ \hline
{ 1} & 1 & 1 & 1 & 1 & 1 & 1 &\\ 
{ 2} & 1 & 1 & 1 & 1 & & &\\ 
{ 3} & 1 & 1 & & & & &\\ 
{ 4} & 1 & 1 & & & & &\\ 
{ 5} & 1 & & & & & &\\ 
{ 8} & 1 & &  &  & & &\\ 
{ {16}} & & & & &  &  &\\ \hline
\end{array} 
$
$
\begin{array}{c|ccccccc|}
{T\rho(\bf 3}) & { 1} & { 2} & { 3} & { 4} & { 5} & { 8} & { {16}}\\ \hline
{ 1} & 1 & 1 & 1 & 1 & 1 & &\\ 
{ 2} & 1 & 1 & & & & &\\ 
{ 3} & 1 & & & & & &\\ 
{ 4} & 1 & & & & & &\\ 
{ 5} & 1 & & & & & &\\ 
{ 8} & & & & & & &\\ 
{ {16}} & & & & & & &\\ \hline
\end{array} 
$
$
\begin{array}{c|ccccccc|}
{T\rho(\bf 4}) & { 1} & { 2} & { 3} & { 4} & { 5} & { 8} & { {16}}\\ \hline
{ 1} & 1 & 1 & 1 & 1 & & &\\ 
{ 2} & 1 & 1 & & & & &\\ 
{ 3} & 1 & & & & & &\\ 
{ 4} & 1 & & & & & &\\ 
{ 5} & & & & & & &\\ 
{ 8} & & & & & & &\\ 
{ {16}} & & & & & & &\\ \hline
\end{array} 
$\\
$
\begin{array}{c|ccccccc|}
{T\rho(\bf 5}) & { 1} & { 2} & { 3} & { 4} & { 5} & { 8} & { {16}}\\ \hline
{ 1} & 1 & 1 & 1 & & & &\\ 
{ 2} & 1 & & & & & &\\ 
{ 3} & 1 & & & & & &\\ 
{ 4} & & & & & & &\\ 
{ 5} & & & & & & &\\ 
{ 8} & & & & & & &\\ 
{ {16}} & & & & & & &\\ \hline
\end{array} 
$
$
\begin{array}{c|ccccccc|}
{T\rho(\bf 8}) & { 1} & { 2} & { 3} & { 4} & { 5} & { 8} & { {16}}\\ \hline
{ 1} & 1 & 1 & & & & &\\ 
{ 2} & 1 & & & & & &\\ 
{ 3} & & & & & & &\\ 
{ 4} & & & & & & &\\ 
{ 5} & & & & & & &\\ 
{ 8} & & & & & & &\\ 
{ {16}} & & & & & & &\\ \hline
\end{array} 
$
$
\begin{array}{c|ccccccc|}
{T\rho(\bf {16}}) & { 1} & { 2} & { 3} & { 4} & { 5} & { 8} & { {16}}\\ \hline
{ 1} & 1 & & & & & &\\ 
{ 2} & & & & & & &\\ 
{ 3} & & & & & & &\\ 
{ 4} & & & & & & &\\ 
{ 5} & & & & & & &\\ 
{ 8} & & & & & & &\\ 
{ {16}} & & & & & & &\\ \hline
\end{array} 
$
\end{footnotesize}
\end{center}

et suivant l'exemple \ref{mertens_exe}, les matrices $ T\rho({\bf u})$ et $T\rho(\bmu)$ :
\begin{center}
\begin{footnotesize}
$
\begin{array}{c|ccccccc|}
 T\rho({\bf u})& { 1} & { 2} & { 3} & { 4} & { 5} & { 8} & { {16}}\\ \hline
 1 & 16& 8& 5& 4& 3& 2&1\\ 
 2 & 8& 4& 2& 2& 1& 1&\\ 
 3 & 5& 2& 1& 1& 1& &\\ 
 4 & 4& 2& 1& 1& & &\\ 
 5 & 3& 1& 1& & & &\\ 
 8 &  2& 1& & & & &\\ 
16 &  1& & & & & &\\ \hline
\end{array} 
$
\hspace{1cm}
$
\begin{array}{c|ccccccc|}
T\rho(\bmu) & { 1} & { 2} & { 3} & { 4} & { 5} & { 8} & { {16}}\\ \hline
 1 & -1& -2& -2& -1& -1& 0&1\\ 
 2 & -2& -1& 0& 0& 1& 1&\\ 
 3 & -2& 0& 1& 1& 1& &\\ 
 4 & -1& 0& 1& 1& & &\\ 
 5 & -1& 1& 1& & & &\\ 
 8 & 0& 1& & & & &\\ 
16 & 1& & & & & &\\ \hline
\end{array} 
$
\end{footnotesize}
\end{center}
\end{exe}

\begin{defi}
\label{UM}
Pour les besoins de la section suivante, nous introduisons les notations :
\begin{equation}
\cU = T\rho({\bf u}) {\ \ \textrm et\ \ } \cM = T\rho({\bmu})
\end{equation}
\end{defi}

\begin{prop}
\label{TuTmu}
Pour $n\in \N^\ast$ fixé, on a :
\begin{equation}
\cU = \left( [n/ij]\right) _{i,j\in \cS} {\ \ \textrm et\ \ }
\cM = \left( M\left( [n/ij]\right) \right) _{i,j\in \cS}
\end{equation}
En particulier, propriété remarquable, la matrice $\cM$ s'obtient par application terme à terme de la fonction de Mertens sur la matrice $\cU$, avec la convention que $M(0)=0$.
\end{prop}
\pr
On a vu dans le corollaire \ref{mertens} que ${\bf u}=\sum_{k\in \mathit{S}}(k-k^-){\bf k}$,
 donc par linéarité :
$$\cU = T\rho({\bf u}) = \sum_{k\in \mathit{S}}(k-k^-)T\rho({\bf k})$$
 et pour $i,j\in\cS$ fixés, 
$$\left( T\rho({\bf u}) \right)_{i,j}=
\sum_{k\in \mathit{S}}(k-k^-)\left( T\rho({\bf k}) \right)_{i,j}$$
D'après la démonstration de la proposition \ref{TM_sym}, on a pour chaque $k\in \cS$ :
$$\left( T\rho({\bf k}) \right)_{ij}=1 \Leftrightarrow  ij\le [n/k]\Leftrightarrow  k\le [n/ij]$$
Il s'ensuit que :
$$\cU_{i,j} = \left( T\rho({\bf u}) \right)_{i,j}=\sum_{k\in \mathit{S},\ k\le [n/ij]}(k-k^-)=[n/ij]$$
ce qu'il fallait démontrer concernant $\cU$.\\
 De même on a (corollaire \ref{mertens}) 
$\bmu=\sum_{k\in \mathit{S}}\left( M(k)-M(k^-)\right) {\bf k}$
et par suite :
 $$\begin{array}{lllll}
\cM & = & T\rho({\bmu}) & = & \sum_{k\in \mathit{S}}\left( M(k)-M(k^-)\right) T\rho({\bf k})\\
\cM_{i,j} & = & \left( T\rho({\bmu}) \right)_{i,j} & = & \sum_{k\in \mathit{S}}\left( M(k)-M(k^-)\right) \left( T\rho({\bf k}) \right)_{i,j}\\
 & & & = &\sum_{k\in \mathit{S},\ k\le [n/ij]}\left( M(k)-M(k^-)\right)  \\
 & & & = & M\left( [n/ij]\right) 
\end{array}$$
\hfill$\square$

\begin{cor}
Pour $n\in \N^\ast$ fixé on a :
\begin{equation}
\vert M(n) \vert \le \Vert \cM_n \Vert
\end{equation}
\end{cor}
\pr
D'après la proposition précédente, on a $M(n) = \cM_{1,1}$ et par ailleurs il est bien connu que pour toute matrice $A$, 
on a $\max_{i,j}\vert a_{i,j}\vert \le \Vert A \Vert$ (voir par exemple \cite{GL} p.57).
\hfill$\square$

\begin{prop}
\label{TUT}
On a 
\begin{equation}
\cM=T\cU^{-1}T
\end{equation}
ce qui signifie que l'on peut construire $\cM$ par inversion de la matrice $\cU$, donc
sans utiliser la propriété remarquable de la proposition \ref{TuTmu}.
\end{prop}
\pr
Rappelons que l'on a $\bmu={\bf u}^{-1}$ (cf. point $\mathit{3.}$ de la proposition \ref{dirichlet} et proposition \ref{bpi_morph}), ce qui entraine d'après la définition de $\cM$ et $\cU$ :
$$\cM=T\rho(\bmu)=T\rho({\bf u})^{-1}=T\left( \cU^{-1}T\right) $$
\hfill$\square$

\begin{rema}
On vérifie facilement les propositions \ref{TuTmu} et \ref{TUT} pour les deux matrices
$\cU=T\rho({\bf u})$ et $\cM=T\rho(\bmu)$ de l'exemple \ref{TuTmu_exe}.
\end{rema}

Dorénavant, lorsque nous voudrons marquer davantage
la dépendance en $n$ des matrices $T,\ \cU,\ \cM$, etc.,
 nous emploierons les notations $T_n,\ \cU_n,\ \cM_n$ à la place de $T,\ \cU,\ \cM$.

\section{Majoration des suites $\Vert T_n \Vert$ et $\VwcM$.}
\label{maj_croiss}

Nous allons nous intéresser, dans cette section et la suivante, à la croissance 
de $\VcM$, $n$ tendant vers l'infini.
Cependant, autant il nous parait difficile d'établir formellement une majoration de $\VcM$ assez fine, c'est-à-dire satisfaisant la conjecture \ref{maj_fondamentale}, autant il est aisé d'obtenir une majoration théorique de la norme de certaines matrices $\wcM$, de conception proche de $\cM$, et que nous allons définir maintenant.
 
\begin{defi}
Pour $n\in\N^\ast$, posons $\wcU = (\wcU_{i,j})_{i,j\in \cS}$ et $\wcM=T\wcU^{-1}T$, où
\begin{equation}
\wcU_{i,j} =
\left\{
\begin{array}{c}
n/ij ,\ \textnormal{si} \ ij\le n \\ 
0 ,\ \textnormal{si} \ ij > n
\end{array} 
\right.
\end{equation}

 Autrement dit $\wcU$ et $\wcM$ sont construites comme $\cU$ et $\cM$ mais sans utilisation de la partie entière $[\ ]$ au dessus de la perdiagonale, (cf. propositions \ref{TuTmu}, \ref{TUT} et définition \ref{UM}).
\end{defi}

\begin{pte}
\label{pteUtilde}
\quad\\
\begin{enumerate}
\item
$T \preccurlyeq \cU \preccurlyeq \wcU$, où le symbole $\preccurlyeq$ entre deux matrices doit se comprendre comme une inégalité terme à terme.
\item
Pour les valeurs élevées de $\cU$, c'est à dire dans le coin nord-ouest des matrices, on a $\cU \simeq \wcU$. \item
Pour les valeurs faibles de $\cU$, c'est à dire sous la perdiagonale et à proximité au dessus, on a $\cU = T$.
\end{enumerate}
\end{pte}

Bien qu'on ne puisse pas tirer des propriétés précédentes de conclusion formelle 
sur le comportement de $\VcM$, il nous parait intéressant d'examiner la croissance des suites $\Vert T_n \Vert$ et $\VwcM$.

\begin{prop}
\label{maj_T}
\begin{equation}
\Vert T_n \Vert = O(\sn)
\end{equation}
\end{prop}

\pr Pour une matrice $A$, posons $\max\vert A \vert = \max_{i,j}\vert a_{i,j}\vert$.
Nous utilisons la majoration $\Vert A \Vert \le s\max\vert A\vert$, 
(voir \cite{GL} p.57), où $s$ désigne la taille de la matrice $A$. Ici, pour $A=T_n$ on a bien sûr d'une part $\max\vert A \vert = 1$ et d'autre part, $s=\#\cS \sim 2\sn$ en vertu de la proposition \ref{diezS}.
\hfill$\square$

\begin{rema}
\label{min_T}

Il est facile de voir que $\lim \inf \Vert T_n \Vert /\sn \ge 1$. Considérons en effet 
le vecteur colonne $w$ de taille $\#\cS$ et constitué de $1$s, et notons $w'$ son transposé.
A cause du fait que $\#\cS \simeq 2\sn$ (proposition \ref{diezS}),
on a d'une part $\Vert w \Vert^2 \simeq 2\sn$ et d'autre part $w'Tw \simeq 2n$. Le rayon spectral de $T_n$, qui est aussi la norme de $T_n$, est donc supérieur à $\dfrac{w'Tw}{\Vert w \Vert^2} \simeq \sn$.
\end{rema}

\begin{lemme}
\label{maj_log}
\begin{equation}
\max\vert \wcM_n \vert=O(\log n)
\end{equation}
\end{lemme}
\pr 
Posons $D=\mathrm{diag}(d_1,\cdots,d_s)$ où $d_k=\sn/k,\ k\in\cS$. On a la forme explicite $\wcU=DTD$, d'où $\wcU^{-1}=D^{-1}T^{-1}D^{-1}$.
Calculons $T^{-1}$, $D^{-1}$ et $\wcU^{-1}$ :
\begin{enumerate}
\item
 $T^{-1}$ est la matrice bi-perdiagonale qui vaut $1$ sur la perdiagonale et $-1$ sur la perdiagonale juste en dessous, (c'est-à-dire décalée d'un cran vers le sud-est).
 \item
$ D^{-1}=(1,\cdots,k,\cdots,n)/\sn=(\ov{n},\cdots,\ov{1})/\sn$, où $k$ parcourt $\cS$ et $\ov{k}=[n/k]$ est le symétrique de $k$ dans la liste $\cS$, (cf. proposition \ref{diezS}).
\item
En parcourant la perdiagonale de $\wcU^{-1}=D^{-1}T^{-1}D^{-1}$, du sud-ouest vers le nord-est, on trouve donc les termes $\dfrac{\ov{k}k}{n}$, $k$ parcourant $\cS$. 
Comme on a $\ov{k}=[n/k]$ il s'ensuit que chacun de ces termes est compris entre $0$ et $1$. Sur la sous-perdiagonale de $\wcU^{-1}$ on trouve les termes $-\dfrac{\ov{k}k^+}{n}$, 
$k$ parcourant $\cS\setminus\{n\}$ et $k^+$ désignant le successeur de $k$ dans $\cS$.
Comme on peut majorer $\dfrac{\ov{k}k^+}{n}$ par $\dfrac{(n/k)k^+}{n} = \dfrac{k^+}{k}\le C$,
$C$ étant la constante obtenue au point $\mathit{3}$ de la proposition \ref{diezS}, chacun de ces termes est donc compris entre $-C$ et $0$.
\end{enumerate}

Enfin, pour en revenir à $\wcM=T\wcU^{-1} T$, on voit qu'obtenir un terme de
$\wcM$ consiste à sommer tous les termes d'une certaine fenêtre rectangulaire de $\wcU^{-1}$. 
Considérant la forme bi-perdiagonale de $\wcU^{-1}$ et le fait que cette matrice est symétrique, on en conclut que chaque coefficient de $\wcM$ est la somme d'au plus :
\begin{itemize}
\item deux sommes latérales, chacune 
de la forme $\dfrac{1}{n}\sum_{k\in\cS,i\le k \le j}\ov{k}(k-k^+)$, $i,j$ fixés dans $\cS$ et $j<\sn$, \item
et d'une somme centrale d'au plus trois termes de $\wcU^{-1}$, chacun compris entre $-C$ et $1$.
\end{itemize}
Lorsque $ k < \sn$, on a vu que $k^+=k+1$, (proposition \ref{diezS}), donc chacune des deux sommes \\
$\dfrac{1}{n}\sum_{k\in\cS,i\le k \le j}\ov{k}(k-k^+)$
se simplifie en $-\dfrac{1}{n}\sum_{k\in\cS,i\le k \le j}\ov{k}$.
En utilisant le point $\mathit{1}$ de la proposition \ref{diezS}, on majore ensuite :
$$\sum_{k\in\cS,i\le k \le j}\ov{k}\le \sum_{1\le k \le \sn}\ov{k}
=\sum_{ 1\le k\le \sn}\left[ n/k\right] \le \sum_{ 1\le k\le \sn} n/k \sim  n\log \sn$$
et on obtient :
$$\max\vert \wcM \vert =0(\log n)$$
\hfill$\square$

\begin{cor}
\label{maj_slog}
\begin{equation}
\VwcM = O(\sn\log n)
\end{equation}
\end{cor}
\pr On utilise le lemme précédent et, comme dans \ref{maj_T},
 l'inégalité $\Vert A \Vert \le s\max\vert A\vert$.
\hfill$\square$

\begin{rema}
Expérimentalement on s'aperçoit que le rapport $\dfrac{\VwcM}{\sn\log n}$ semble tendre vers 0,
bien qu'il soit assez facile de prouver que $\dfrac{\max\vert \wcM_n \vert}{\log n}$ a une limite strictement positive (résultat donc plus précis que dans le lemme \ref{maj_log}). Ceci indiquerait donc que la majoration $\Vert A \Vert \le s\max\vert A\vert$ est un peu trop sévère dans notre cas.
\end{rema}

\section{Expérimentation.}
Dans cette section nous montrons le résultat de quelques expérimentations numériques
sur la croissance comparée des suites $\VcM$, $M(n)$, $\VwcM$ et $\Vert T_n \Vert$ lorsque $n$ croit vers l'infini.

\subsection{Croissance comparée de $\VcM$ et $M(n)$.}
 Le graphique ci-dessous montre les courbes de $\VcM/\sn$ en noir, et $M(n)/\sn$ en bleu, pour les valeurs de $n$ s'échelonnant de $5000$ en $5000$ jusqu'à $10^6$.

\begin{figure}[!ht]
\label{num_results_41}
\centering
\includegraphics[width=1.0\textwidth, height=0.6\textwidth]{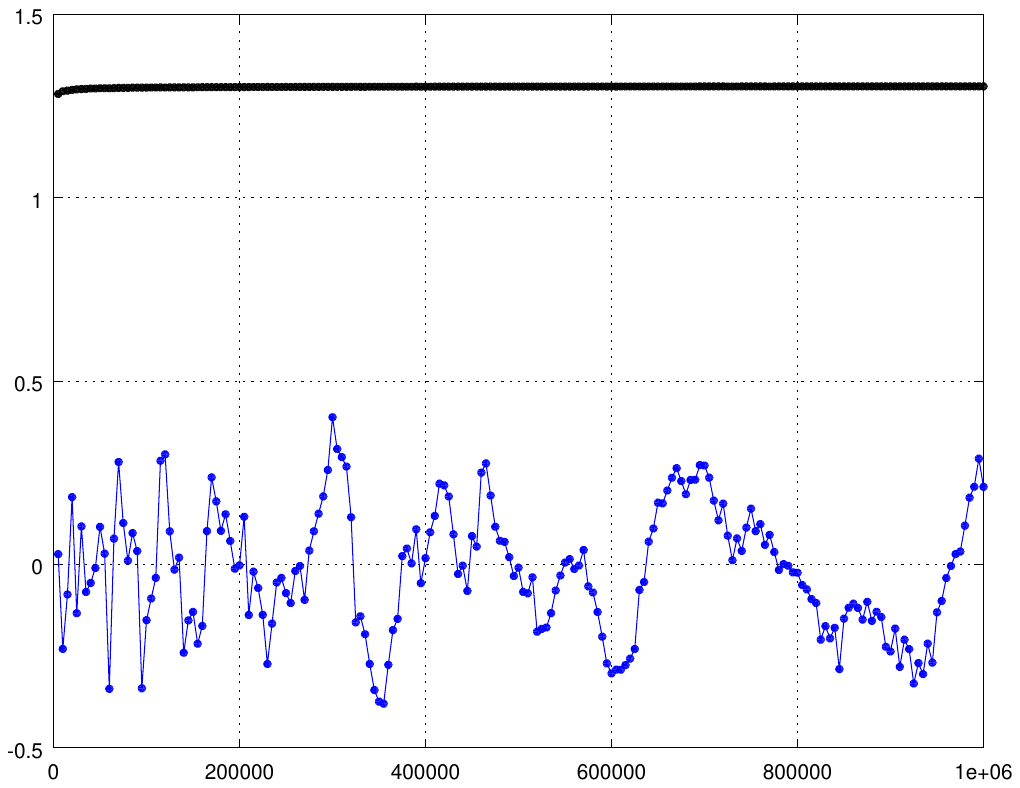}
\end{figure} 

On constate une grande régularité de $\VcM/\sn$, contrastant avec le caractère chaotique de  $M(n)/\sn$. Autre observation importante : $\VcM/\sn$ semble croître assez lentement, ce que nous allons tenter de préciser dans les graphiques suivants.

\subsection{Croissance comparée de $\VcM$, $\VwcM$ et $\VT$.}
\label{croissance_MMT}
 Nous comparons maintenant la croissance de $\VcM$ avec les suites $\VwcM$ et $\VT$ rencontrées dans la section précédente. En ce qui concerne la croissance de ces deux suites, rappelons que l'on a $\VwcM = O(\sn\log n)$ et $\Vert T_n \Vert = O(\sn)$, (cf. corollaire \ref{maj_slog} et proposition \ref{maj_T}).
 La figure de gauche ci-dessous montre les courbes respectives de $\VcM /\sn$ en noir,
 de $\VwcM /\sn$ en rouge, et de $\VT /\sn$ en bleu.
A droite nous montrons les mêmes suites mais normalisées chacune par leur valeur pour $n = 5 \times 10^5$, ce qui permet de mieux comparer leurs croissances.
\begin{figure}[!ht]
\label{num_results_42}
\centering
\includegraphics[width=1.0\textwidth, height=0.6\textwidth]{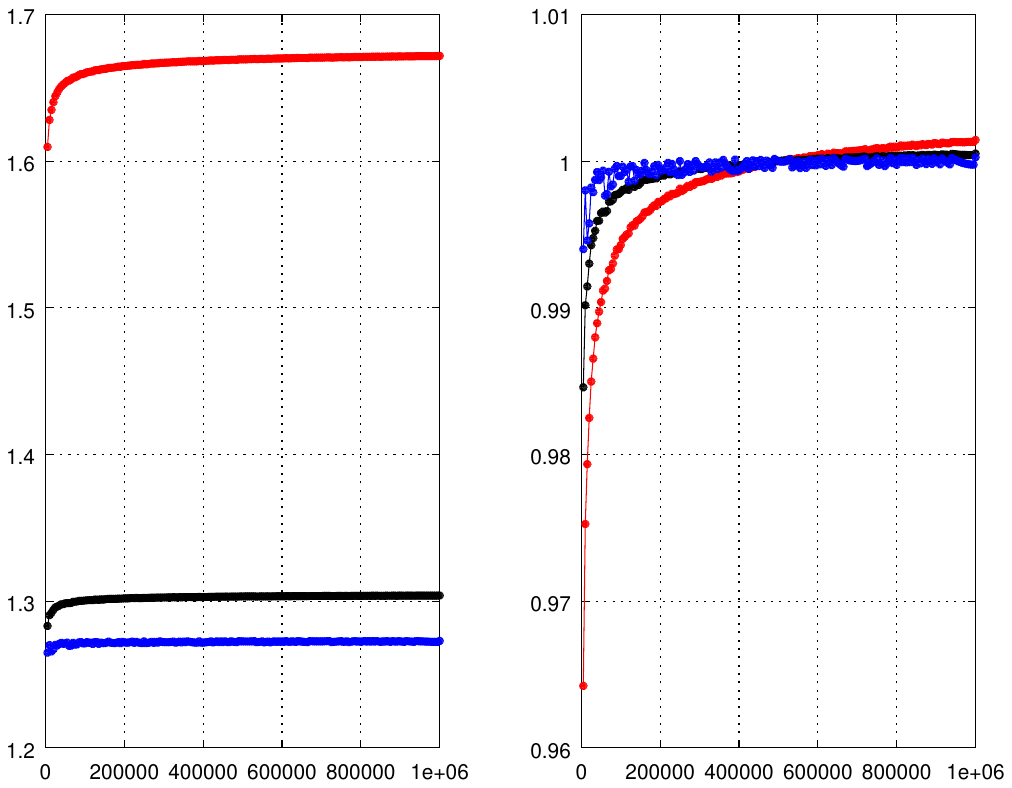}
\end{figure}

On obsverve que $\VcM$ semble croître plus vite que $\VT$ mais moins vite que $\VwcM$. Comme nous savons que la croissance de $\VT$ est comparable à celle de $\sn$ (cf. proposition \ref{maj_T} et remarque \ref{min_T}), nous allons maintenant examiner graphiquement les suites $\VcM/\sn$ et $\VcM/\VwcM$.

\pagebreak

\subsection{Courbes de $\VcM/\sn$ et $\VcM/\VwcM$.}
Le graphique ci-dessous montre, à gauche la courbe de $\VcM / \sn$, et à droite celle de $\VcM / {\VwcM}$.

\begin{figure}[!ht]
\label{num_results_43}
\centering
\includegraphics[width=1.0\textwidth, height=0.6\textwidth]{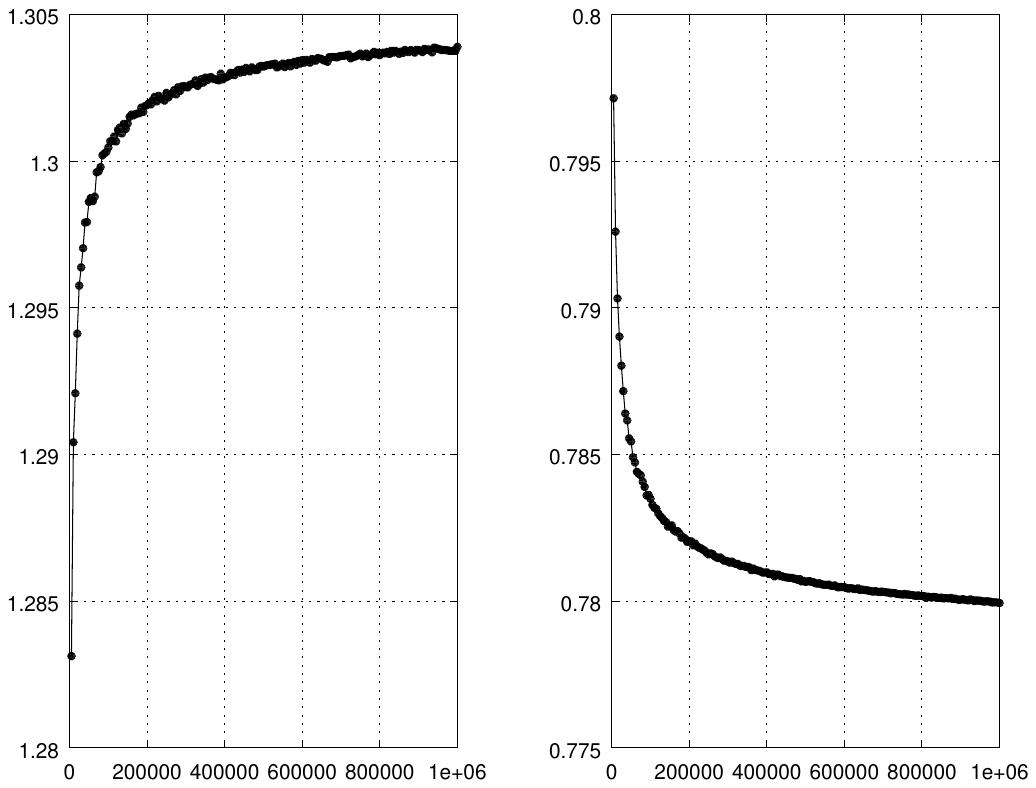}
\end{figure}

Observation et discussion :
\begin{enumerate}
\item
A l'observation de la courbe de gauche, sans envisager toutes les possibilités pour \\
$\lim\inf \VcM/\sn$ et $\lim\sup {\VcM/\sn}$, nous retenons quelques hypothèses :
\begin{itemize}
\item  $\lim_{n\rightarrow\infty}\ \VcM/\sn = +\infty$, qui signifie que $\VcM$ croît plus vite que $\sn$. Ceci serait en accord avec les conjectures avancées par différents auteurs sur la croissance de $\vert M(n) \vert$. Par exemple T. Kotnik et J. van de Lune \cite{KL2} conjecturent que l'on a :
$$\lim\sup_{n\rightarrow\infty}\vert M(n) \vert / \sqrt{n\log\log\log n} = C$$
 où $C$ est une  constante non nulle. On trouvera dans \cite{KL2} une discussion et des références sur des estimations de $\vert M(n) \vert$ par d'autres auteurs.
\item
$ \lim_{n\rightarrow\infty}\ \VcM/\sn  < +\infty$, qui signifie que les croissances de 
$\VcM$ et $\sn$ sont comparables. Cette hypothèse contredirait les conjectures que l'on vient d'évoquer dans le point précédent, mais ne peut être formellement exclue.
\item Enfin remarquons que l'hypothèse $\lim_{n\rightarrow\infty} \VcM/\sn = 0$, qui signifie que $\VcM$ croît moins vite que $\VT$, est exclue car elle contredirait le fait établi dans \cite{OR} que \\
$\lim\sup_{n\rightarrow\infty} \vert M(n) \vert /\sn > 1.06$. Cette hypothèse serait d'ailleurs peu en accord avec les tendances apparaissant dans le graphique.
\end{itemize}

\item En ce qui concerne la courbe de droite, nous envisageons principalement deux hypothèses :
\begin{itemize}
\item $\lim_{n\rightarrow\infty}\ \VcM / \VwcM < +\infty$. D'après la courbe, cette hypothèse semble la plus probable, et en tenant compte du corollaire \ref{maj_slog}, elle implique la conjecture \ref{maj_fondamentale}.
\item L'hypothèse $\lim_{n\rightarrow\infty}\ \VcM / \VwcM = +\infty$, qui signifie que $\VcM$ croît plus vite que $\VwcM$, semble peu probable d'après la courbe et les observations de la section 
\ref{croissance_MMT}, et nous la rejetons, bien qu'elle ne contredise pas définitivement la conjecture \ref{maj_fondamentale}.
\end{itemize}
\end{enumerate}

\section{Conclusion.}

Nous avons construit une suite de matrices $\cM$ vérifiant $\forall n\in\N^\ast \vert M(n) \vert \le \VcM$
où $M$ désigne la fonction de Mertens sur les entiers. Une expérimentation numérique sur l'intervalle $n\in\left[ 5\times 10^3,  10^6\right]$ nous a incité à énoncer la conjecture \ref{maj_fondamentale} : 
\begin{equation}
 \forall \epsilon > 0, \VcM =O(n^{1/2+\epsilon})
\end{equation}
qui implique la conjecture de Riemann. A aucun moment nous n'avons employé de variable complexe ni fait usage des techniques habituelles de la théorie analytique des nombres. Bien que quelques auteurs aient déjà envisagé une approche matricielle à ce problème (voir par exemple \cite{BJ}), le travail présenté ici nous amène à penser que les méthodes spectrales en analyse matricielle pourraient jouer un rôle accru dans la recherche sur l'hypothèse de Riemann.

\bibliographystyle{abbrv}

\hspace{-0.65cm} Dpt. de Mathématiques, Université Paris 13, 99 av J.B. Clément, 93430 Villetaneuse, France\\
Email adress : cardinal@math.univ-paris13.fr

\end{document}